\documentclass{article}
\usepackage{amsmath,amssymb,amsthm}

\begin{document}

\title{ Linear Algebra Estimates}         
\author{Hatice \c Sahino\u glu}        
\date{\today}          
\maketitle

\newtheorem{thm}{Theorem}[section]
\newtheorem{cor}[thm]{Corollary}
\newtheorem{lem}[thm]{Lemma}
\newtheorem{prop}[thm]{Proposition}
\newtheorem{definition}[thm]{Definition}
\theoremstyle{definition}
\newtheorem{rmk}[thm]{Remark}

\newcommand{\Aut}{\operatorname{Aut}}

\begin{abstract}
In this paper we give a generalization of a linear algebra estimate that occurs in the paper \cite{RS}, by Michael Rosen and Joseph H. Silverman. In \cite{RS} authors give a bound for the size of a submodule of $(\mathbb{Z}/n \mathbb{Z})^2$ in terms of a power of the index of any subgroup of automorphism group of $(\mathbb{Z}/n \mathbb{Z})^2$  which is acting in an abelian way on that submodule, meaning that given $A$ and $B$ as any two elements in the automorphism group $AB-BA$ annihilates all elements of the submodule. We will give a similar estimate for finite submodules of arbitrary dimension $m$ and subgroups of general linear group acting on them. Later we will derive the analog of this result for the case of subgroups of the symplectic group acting on finite submodules in an abelian fashion.
\end{abstract}
\section{Introduction}
Given a not necessarily commutative group acting on a module, if we force the group to act in an abelian way on this module we should either take a small module or restrict the action group to a subgroup of it as we do not expect a big subgroup of a highly non-commutative group to act in an abelian way on a submodule of large size. The reason for that is if the submodule is small enough the group is more likely to be able commute while acting on this module. In the same way if we get rid of some elements, those that do not commute with the other elements while acting  on the module, of the matrix group the remaining elements will act in an abelian way. Hence we give a relation on the submodule in terms of the index of the subgroup of a matrix group that is acting on the submodule in an abelian way for two different matrix groups.  One of these is the general linear group in arbitrary dimension and the other one is the symplectic group in general dimension.

Our first goal is to prove the following proposition.
\begin{prop}\label{pro1}
Suppose
\begin{itemize}
\item{}$n,m$ are positive integers greater than $1$.
\item{}$V$ is a free $(\mathbb{Z}/n \mathbb{Z})$- module of rank $m$.
\item{}$\Gamma$ is a subgroup of $Aut(V)\cong GL_m(\mathbb{Z}/n \mathbb{Z})$.
\item $W$ is a $\Gamma$-invariant $\mathbb{Z}/n \mathbb{Z}$- submodule of $V$.
And also that the action of $\Gamma$ on $W$ is abelian. That's to say $\Gamma|_W $ is an abelian subgroup of $Aut(W)$.
\item{}
Let
\begin {equation}
 I(\Gamma):=(Aut(V): \Gamma)=\text{the index of } \Gamma \text { in } Aut(V)
\end{equation}
\end{itemize}
Then we have
\begin{equation}
\|W\|\leq I(\Gamma)^{1+(\frac{m}{2})}\zeta(2)^{(m-1)(1+\frac{m}{2})}
\end{equation}
where $\|W\|$ is the magnitude of the set $W$.
\end{prop}
\cite{RS}  also gives a lower bound for the magnitude of the $\|W\|$ in terms of the magnitude of $I(\Gamma)$. For dimension two, it is shown there exists $W$ such that
$$\|W\|\geq I(\Gamma)^{\frac{4}{3}}.$$

For higher dimension $m$, using the same type of submodules with \cite{RS} one can show there exist $W$ such that
\begin{equation}
\|W\|\geq I(\Gamma)^{1+(\frac{1}{m+1})}.
\end{equation}
We will not go into the details since it is a straight forward deduction.
\\
After we prove the Proposition (\ref{pro1}) above, for the purpose of generalizing the  work of \cite{RS} to higher dimensional abelian varieties we will need a version of Proposition (\ref{pro1}) for symplectic matrices. Hence we will show that the following theorem is also true.
\begin{thm} 
Given a submodule $W$ of $V$ where $V$ is the free $Z/nZ$ submodule of dimension $2m$, where $m$ is an integer, let  $\Gamma$ be a subgroup of the symplectic group $Sp(2m, Z/nZ)$ that is acting in an abelian way on $W$ and let  $I(\Gamma)$ be the index of $\Gamma$ in the symplectic group, 
then 
$$|W| \leq I(\Gamma)^{(2+2m)}\zeta(2)^{(2+2m)}.$$  
\end{thm}

\section{General Dimension Case }       

Let $W$ be a $\mathbb{Z}/n \mathbb{Z}$ submodule of $(\mathbb{Z}/n \mathbb{Z})^m$
and let $\Gamma$ be a subgroup of $GL_m(\mathbb{Z}/n \mathbb{Z})$ that acts in an abelian way on $W$.
Using the notation introduced earlier, let $\|W\|$ be the order of $W$ and
$$I(\Gamma)=(Aut(V): \Gamma)=\text{ the index of } \Gamma \text { in } Aut(V).$$ 
We will prove that,

\begin{equation}\label {main}
\|W\|\leq I(\Gamma)^{1+(\frac{m}{2})}\zeta(2)^{(m-1)(1+\frac{m}{2})}.
\end{equation}

\begin{proof}

Let's decompose $W$ and $\Gamma$ into $l$-primary parts. $V=\oplus_{l|n}V_l$ and $W=\oplus_{l|n}W_l$ where $V_l=V\otimes \mathbb{Z}_l$ and similarly $W_l=W\otimes \mathbb{Z}_l$.
\label{decomposition} $Aut(V)=\displaystyle \oplus _{l|n}Aut(V_l)$ and $\Gamma=\displaystyle \oplus_{l|n}\Gamma_l$ where $\Gamma_l=Image(\Gamma \rightarrow Aut(V_l))$.

If we can show that the following relation holds for the $l$-primary parts of $W$ and $\Gamma$, i.e,
\begin{equation}\label {mainsub}
\|W_l\|\leq I(\Gamma_l)^{1+(\frac{m}{2})}c_l.
\end{equation}
where $\prod c_l$ converges, then we get

\begin{align}
\|W\|&\leq \prod _{l|n}\|W_l\| \leq \prod_{l|n}I(\Gamma_l)^{1+(\frac{m}{2})}c_l =I(\Gamma)^{1+(\frac{m}{2})}c.
\end{align}

Therefore for the rest of the proof we will assume that $n$ is a power of a prime number, say $n=l^e$ and try to show that $\|W\| \leq I(\Gamma)^{1+(\frac{m}{2})} {c_l}$.  
By choosing an appropriate basis for $W$ let 
$$W={l^i\mathbb{Z}}/{l^e\mathbb{Z}}\times {l^{i+j_1}\mathbb{Z}}/{l^e\mathbb{Z}}\times {l^{i+j_1+j_2}\mathbb{Z}}/{l^e\mathbb{Z}}\times...{l^{i+j_1+j_2+...j_{m-1}}\mathbb{Z}}/{l^e\mathbb{Z}}.$$

Given any matrix $A$ in $\Gamma$ as below,

$$A= \left( \begin{array}{cccccc}
a_{11} & a_{12} & a_{13}& . & . & a_{1m} \\
a_{21} & a_{22} & a_{23}& . & . & a_{2m} \\
a_{31} & a_{32} & a_{33}& . & . & a_{3m} \\
. &.&.&.&.&.\\
. &.&.&.&.&.\\
a_{m1} & a_{m2} & a_{m3}& . & . & a_{mm}\end{array} \right).$$ 
The condition $\Gamma W=W$ implies that  multiplying first row with the first basis element we should get a number congruent to $0$ modulo $l^{i}$; second row multiplied with the first basis element should give a number congruent to $0$ modulo $l^{i+j_1}$ and $n$'th row multiplied with the first basis element should give a number congruent to $0$ modulo $l^{i+j_1+...+j_{n-1}}$. In the same manner $n$'th row of $A$ multiplied with the $k$'th basis element should give a number congruent to $0$ modulo $l^{i+j_1+...+j_{n-1}}$. 
Let us denote $l^{j_a+j_{a+1}+...j_b}$ as $l^{j_{(a)(b)}}$ for sake of simplicity in notation.
Therefore $A$ should be in the following form:
$$ \left( \begin{array}{cccccc}
a_{11} & a_{12} & a_{13}& . & . & a_{1m}\\
l^{j_1}\bar{a}_{21} & a_{22} & a_{23}& . & . &a_{2m} \\
l^{j_{(1)(2)}}\bar{a}_{31}& l^{j_2}\bar{a}_{32}&a_{33} & .& . &a_{3m} \\
. &.&.&.&.&.\\
. &.&.&.&.&.\\
l^{j_{(1)(m-1)}}\bar{a}_{m1} &l^{j_{(2)(m-1)}} \bar{a}_{m2} & l^{j_{(3)(m-1)}} \bar{a}_{m3}& .&.&a_{mm} \end{array} \right).$$

Let $$J=\min \left \{ ord_l (a_{21}):  \left( \begin{array}{cccccc}
a_{11} & a_{12} & a_{13}& . & . & a_{1m}\\
a_{21} & a_{22} & a_{23}& . & . &a_{2m} \\
a_{31}& a_{32}&a_{33} & .& . &a_{3m} \\
. &.&.&.&.&.\\
. &.&.&.&.&.\\
a_{m1} &a_{m2} & a_{m3}& .&.&a_{mm} \end{array}\right ) \in \Gamma \right \}.$$
And if $\Gamma$ consists of matrices all with $a_{21}=0$ we set $J=e$.
Hence there exists a matrix $A$ in $\Gamma$ with $a_{21}=l^J \bar{a}_{21}$ such that $\bar{a}_{21}\not\equiv 0 \pmod {l}$. Let $A\in \Gamma$ be
$$ \left( \begin{array}{cccccc}
a_{11} & a_{12} & a_{13}& . & . & a_{1m}\\
l^{J}\bar{a}_{21} & a_{22} & a_{23}& . & . &a_{2m} \\
a_{31}& a_{32}&a_{33} & .& . &a_{3m} \\
. &.&.&.&.&.\\
. &.&.&.&.&.\\
a_{m1} &a_{m2} & a_{m3}& .&.&a_{mm} \end{array} \right) $$
with $\bar{a}_{21}\not\equiv 0 \pmod {l}$.
\\
For any $B$ in $\Gamma$ given as
$$ \left( \begin{array}{cccccc}
a'_{11} & a'_{12} & a'_{13}& . & . & a'_{1m}\\
l^{J}\bar{a}'_{21} & a'_{22} & a'_{23}& . & . &a'_{2m} \\
a'_{31}& a'_{32}&a'_{33} & .& . &a'_{3m} \\
. &.&.&.&.&.\\
. &.&.&.&.&.\\
a'_{m1} &a'_{m2} & a'_{m3}& .&.&a'_{mm}
 \end{array} \right) $$

$B$ commutes with $A$ on $W$, i.e

\begin{equation} \label{Ma}
 (AB-BA)  \left( \begin{array}{cccccc}
l^i& 0 & 0& . & . & 0 \\
0 & l^{i+j_1} & 0& . & . &0 \\
0& 0 &l^{i+j_{(1)(2)}} & .& . &0 \\
. &.&.&.&.&.\\
. &.&.&.&.&.\\
0 & 0 & 0& .&.&l^{i+j_{(1)(m-1)}}\end{array}\right) \equiv 0 \pmod {l^e}. 
\end{equation} 
We will replace $\bar{a}_{21}$ with $a_{21}$ and $\bar{a}'_{21}$ with $a'_{21}$ in the representation of $A$ and $B$ for sake of simplicity in notation.\\

Without doing the whole matrix multiplication, we will pass to the counting of $B$'s in $\Gamma$. 
Given $A$ and any $B$ as above there are matrices $X$ such that the following congruence holds. 
\begin{equation}\label{forXn}
l^{i+J}[(a'_{11}a_{21}-a_{11}a'_{21})I+a'_{21}A]\equiv l^{i+J}[a_{21}B+X]\,\pmod { l^e}
\end{equation}
$X$ should be as follows,
$$X \equiv  \left( \begin{array}{cccccc}
0&a'_{21}a_{12}-a_{21}a'_{12}&a'_{21}a_{13}-a_{21} a'_{13}& * & * & a'_{21}a_{1n}-a_{21}a'_{1n} \\
0 &a'_{21}a_{11}-a_{21}a'_{11}+a'_{21}a_{22}-a_{21}a'_{22}  & *& * & * &a'_{21}a_{2n}-a_{21}a'_{2n} \\
*&(a'_{21}a_{32}-a_{21}a'_{32}) & *& * &*&* \\
* &*&*&*&*&*\\
* &*&*&*&*&*\\
* &(a'_{21}a_{m2}-a_{21}a'_{m2}) & *& *&*&a'_{21} a_{mm}- a_{21} a'_{mm}\end{array}\right) \equiv 0 \pmod{ l^{e-i-J}}. $$  
The goal is to write entries with $a'_{21},a'_{22},a'_{32},...,a'_{m2}$ that occurs in $X$ in terms of the other entries of $X$ by using (\ref{Ma}). 
\\
First row of $A$ multiplied with the first column of $B$ gives: 
$$a_{11}a'_{11}+ l^{J}(a_{12}a'_{21})+(a_{13}a'_{31})+....+(a_{1m} a'_{m1}).$$ 
After taking the symmetric difference $AB-BA$ and multiplying with the basis matrix of $W$ one gets,
$$
\begin{aligned}
&l^{i+J}(a_{12}a'_{21}-a'_{12}a_{21})+l^{i} (a_{31}a'_{31}-a'_{31}a_{31})+...+l^{i}(a_{1m} a'_{m1}-a'_{1m} a_{m1}) 
\end{aligned}
.$$
This symmetric difference is congruent to $0 \pmod{l^e}$. This gives that $l^{i+J}(a'_{21}a_{12}-a_{21}a'_{12})$ can be written in terms of entries $a_{ij}$ and $a'_{ij}$, for $j\neq 2$, $\pmod{l^{e}}$ where $a_{ij}$ and $a'_{ij}$ are entries of matrices of $A$ and $B$ respectively.

For the second entry of the second row of $X$, consider second row of $(AB-BA)$ multiplied with the first basis vector. The first entry of the second row of $(AB-BA)$ is the following symmetric difference 
$$
\begin{aligned}
&l^{J}(a_{21}a'_{11}-a'_{21}a_{11})+l^{J}(a_{22}a'_{21}-a'_{22}a_{21})+(a_{23}a'_{31}-a'_{23}a_{31})+...+(a_{2m} a'_{m1}-a'_{2m}a_{m1}).
\end{aligned}
$$ 
Multiplying this entry with the first basis element we get that 
$$(a'_{21}a_{11}-a_{21}a'_{11})+(a'_{21}a_{22}-a_{21}a'_{22})$$

can be written in terms of entries $a_{ij}$ and $a'_{ij}$ for $i,j>2$, in $\bmod{ l^{e-i-J}}$ .

For $k$'th entry in the second column of matrix $X$, consider the first entry of the $k$th row of the matrix $(AB-BA)$. This entry is 
$$
\begin{aligned}
&(a_{k1} a'_{11}-a'_{k1} a_{11})+l^{J}(a_{k2}a'_{21}-a'_{k2}a_{21})+\\
&(a_{k3}a'_{31}-a'_{k3}a_{31})+...(a_{km}a'_{m1}-a'_{km}a_{m1}).
\end{aligned}
$$
From (\ref{Ma}) this identity should be congruent to $0$ modulo ${l^{e-i}}$. Hence $a'_{21}a_{k2}-a_{21}a'_{k2}$ can be written in terms of entries in $k$'th row and the first column modulo  $l^{e-i-J}$.
Therefore terms with $a'_{12}, a'_{22}, a'_{32},...,a'_{m2}$ in the representation of $X$ are determined by the remaining entries of $X$ modulo $l^{e-i-J}$. 
Remember that in equation \ref{forXn} we had
\begin{equation*}
l^{i+J}[(a'_{11}a_{21}-a_{11}a'_{21})I+a'_{21}A]\equiv l^{i+J}[a_{21}B+X] \; \pmod {l^e}.
\end{equation*}
Since $a_{21}\neq 0 \pmod{l}$, dividing the equation above by $a_{21}$ we get:
\begin{equation}
B \equiv xI+yA-X \, \text{mod}\, l^{e-i-J}\,\text{ for } x,y\in {\mathbb{Z}}/{l^{e-i-J}\mathbb{Z}} \, \\
\text{ and }\, X \, \text {is as above}.
\end{equation}
This is equivalent to the following:
 \begin{equation} 
B=xI+yA-X+l^{e-i-J}Z \text{  where  } x,y\in \mathbb{Z}/{l^{e-i-J}\mathbb{Z}} \text{ and } Z\in M_m(\mathbb{Z}/{l^{i+J}\mathbb{Z}}).
\end{equation}
Note that choices of $x,y$ are determined by the choice of $a'_{11},a'_{21}.$
\\
Now we have enough information to count number of $B$'s and then get an upper bound for $\Gamma$. 
\begin{align}
\nonumber |\Gamma|&\leq (\text{\# of choices of }(x,y))({\text{\# of choices of }X}) (\text{\# of choices of }Z)\\
&=l^{2(e-i-J)}l^{(m^2-m-2)(e-i-J)}l^{(m^2)(i+J)}\\
&=l^{(m^2-m)(e-i-J)+m^{2}(i+J)}.
\end{align}
implies
\begin{align} \label{indexm1}
I(\Gamma)= |\text{Aut}(GL_m(\mathbb{Z}/{l^{e}\mathbb{Z}})/\Gamma| &\geq \frac{ l^{m^{2}e}(1-l^{-1})(1-l^{-2})...(1-l^{-m})}{l^{m^2e-m(e-i-J)}}\nonumber \\
&=l^{m(e-i-J)}(1-l^{-1})(1-l^{-2})...(1-l^{-m})
\end{align}

On the other hand, elements of B's in $\Gamma$ are in the following form due to divisibility conditions  that come from the action of $\Gamma$ on $W$ being invariant.

$$ \left( \begin{array}{cccccc}
a_{11} & a_{12} & a_{13}& . & . & a_{1m}\\
l^{j_1}\bar{a}_{21} & a_{22} & a_{23}& . & . &a_{2m} \\
l^{j_{(1)(2)}}\bar{a}_{31}& l^{j_2}\bar{a}_{32}&a_{33} & .& . &a_{3m} \\
. &.&.&.&.&.\\
. &.&.&.&.&.\\
l^{j_{(1)(m-1)}}\bar{a}_{m1} &l^{j_{(2)(m-1)}} \bar{a}_{m2} & l^{j_{(3)(m-1)}} \bar{a}_{m3}& .&.&a_{mm} \end{array} \right).$$
Note that $j_{(a)(b)}$ is the short cut for the sum $j_a+j_{a+1}+...j_b$.

We will get rid of bar sign for the convenience of notation.

From this we get the following inequality for the index of $\Gamma$ in the general linear group.
\begin{align}\label{indexmsub2}
I(\Gamma) &\geq \frac{l^{m^{2}e}(1-l^{-1})(1-l^{-2})...(1-l^{-m})}{l^{m^2e-J-(m-2)j_1-2(m-2)j_2-3(m-3)j_3...}(1-l^{-1})^m}\nonumber
\\&=l^{J+(m-2)j_1+2(m-2)j_2+....}(1-l^{-1})^{(1-m)}(1-l^{-2}) (1-l^{-3})...(1-l^{-m})
\end{align}
Now we will investigate how large $J$ compared to  $(j_1+j_2+...j_{m-1})$ is. $\Gamma$ is closed under multiplication. For any matrix $A$ in $\Gamma$ given as below; 

\[ \left( \begin{array}{cccccc}
a_{11} & a_{12} & a_{13}& . & . & a_{1m}\\
l^{J}a_{21} & a_{22} & a_{23}& . & . &a_{2m} \\
l^{j_{(1)(2)}}a_{31}& l^{j_2}a_{32}&a_{33} & .& . &a_{3m} \\
. &.&.&.&.&.\\
. &.&.&.&.&.\\
l^{j_{(1)(m-1)}} a_{m1} &l^{j_{(2)(m-1)}} a_{m2} & l^{j_{(3)(m-1)}} a_{m3}& .&.&a_{mm}
 \end{array} \right)\]
$A^2\in \Gamma$.  Looking at the second row first entry of $A^2$ we should see a number divisible by $l^{J}$. This number is 
\begin{equation}\label{above}
l^{J}(a_{11}a_{21}+a_{21}a_{22})+l^{j_1+j_2}a_{31}a_{23}+l^{j_1+j_2+j_3}a_{41}a_{24}+....+l^{j_1+j_2+j_3+....j_{m-1}}a_{m1}a_{2m}.
\end{equation} 
Assume that $J\geq (j_1+j_2+...j_{m-1})$. Since the quantity (\ref{above}) should also be congruent to $0$ modulo $l^{j_1+j_2+j_3}$ we conclude that $a_{31}a_{23}$ divisible by $l^{j_3}$. Similarly it should be congruent to $0$ modulo $l^{j_1+j_2+j_3+j_4}$. This implies that 
$$l^{j_4}| (a_{31}a_{23}/l^{j_3}+a_{41}a_{24}).$$ 
Continuing in the same way we finally get that 
$$l^{J-(j_1+j_2+...j_{m-1})}|(\frac{1}{l^{j_{m-1}}}(...\frac{1}{l^{j_4}}(\frac{1}{l^{j_3}} a_{31}a_{23}+a_{41}a_{24})+....)+a_{n1}a_{2n}).$$
This last restriction tightens the set of ${a_{31},a_{41}...a_{n1},a_{23}...a_{2n} }$ by $l^{J-(j_1+j_2+...j_{m-1})}$. To see this note that for every choice of ${a_{31},a_{41}...a_{n1},a_{23}...a_{2{n-1}} }$; $a_{2n}$ is determined up to a multiple of $l^{J-(j_1+j_2+...j_{m-1})}$. Hence the assumption $(j_1+j_2+...j_{m-1})\leq J$ leads to an extra factor of $l^{J-(j_1+j_2+...j_{m-1})}$ in the index $I(\Gamma)$ above. Therefore the index in (\ref{indexmsub2}) can be bounded as follows.
\begin{align}\label{indexm2}
I(\Gamma)\geq l^{2J+(m-2)j_1+(m-3)j_2+....}(1-l^{-1})^{(1-m)}(1-l^{-2}) (1-l^{-3})...(1-l^{-m}).
\end{align}

Finally multiplying the first power of (\ref{indexm1}) and $m/2$ 'th power of (\ref{indexm2}),
\begin{align} \label{last}
I(\Gamma)^{m/2+1} &\geq l^{me-mi+c} (1-l^{-1})^{1+\frac{m}{2}(1-m)}[(1-l^{-2}) (1-l^{-3})...(1-l^{-m})]^{1+m/2} \text{ where } c \geq 0.
\end{align}
 
Noting that
$$\|W\|=l^{me-mi-(m-1)j_1-c_2}\text{ where }c_2 \geq 0$$
(\ref{last}) implies that the following relation is true
\begin{equation}
I(\Gamma)^{1+\frac{m}{2}}\geq \|W\|(1-l^{-1})^{1+\frac{m}{2}(1-m)}[(1-l^{-2}) (1-l^{-3})...(1-l^{-m})]^{(1+m/2)}.
\end{equation}
Equivalently,
\begin{align}
&I(\Gamma)^{1+\frac{m}{2}}\geq \|W\|[(1-l^{-2}) (1-l^{-3})...(1-l^{-m})]^{1+m/2}\text{ implies}\\
&I(\Gamma)^{1+\frac{m}{2}}\geq \|W\|(1-l^{-2})^{(m-1)(1+\frac{m}{2})}.
\end{align}
Noting that
$$\prod_{l}[(1-l^{-2}) (1-l^{-3})...(1-l^{-m})]^{1+m/2} \geq \zeta(2)^{-(m-1)(1+\frac{m}{2})},$$
if we drop the constraint on $n$ being a prime power we have 
\begin{equation}
\|W\| \leq I(\Gamma)^{1+\frac{m}{2}}\zeta(2)^{(m-1)(1+\frac{m}{2})}.
\end{equation}
This finishes the proof.

 \end{proof}

\section{Symplectic Group Case}
\begin{definition}Symplectic group is the set of  $2m \times 2m$ matrices $M$ such that $M^{T}\Omega M =\Omega$ where $\Omega$ is the block matrix
$$
 \left ( \begin{array}{cc}
0 & I_m  \\
-I_m & 0 \end{array} \right)
$$ with $I_m$, the $m\times m$ identity matrix.
\end{definition}
Let $M$ be a  matrix in the symplectic group with entries 
$$\left ( \begin{array}{cc}
A & B  \\
C & D \end{array} \right)
$$ where $A,B,C,D$ are $m\times m$ matrices.
The symplectic property imposes that 
\begin{itemize}
\item{} $A^{T} D-C^{T}B=I$,
\item{} $A^{T} C=C^{T}A$, 
\item{} $D^{T} B=B^{T}D$.
\end{itemize}

The property $A^T D-C^T B= I$ gives $m^2$ relations on the elements of the matrices $A,B,C,D$.
The property $A^T C=C^T A$ gives $\frac{m^2-m}{2}$ relations on the entries of $A,C$.
Similarly the property $D^T B=B^T D$ brings $\frac{m^2-m}{2}$ relations on the entries of $A,C$. 
Now let's assume that we are working with symplectic matrices of dimension $2m$ on the ring $Z/nZ$. This group is denoted as $Sp(2m, Z/nZ)$.
\begin{thm} 
Given a submodule $W$ of $V$ where $V$ is the free $Z/nZ$ submodule of dimension $2m$, let  $\Gamma$ be a subgroup of the symplectic group $Sp(2m, Z/nZ)$ that is acting in an abelian way on $W$ and let  $I(\Gamma)$ be the index of $\Gamma$ in the symplectic group, 
then 
$$|W|\leq I(\Gamma)^{(2+2m)}\zeta(2)^{(2+2m)}.$$  
\end{thm}
In the previous section where we worked with $GL(2m,Z/nZ)$, we could use the commutativity property to derive some constraints on the entries of matrices. In the case of symplectic matrices one should be careful since the commutativity property may not bring a totally new set of relations.

Instead we start with a column with a decent number of free entries (i.e a column whose entries take all possible range of values as the matrices vary in the whole symplectic group) and on these free elements we can use the implications of the commutativity to derive extra constraints.  
The group of symplectic matrices $Sp(2m,k)$ has dimension $2m^2+m$ (Look at \cite{Dav}) over the field $k$. This is also true when $k$ is a ring (Look at the Appendix). By pigeon hole principle there should be at least one column, different than the first one, with at least $m$ free entries, say this is the $k$'th column. We will show that some of these free entries will be eliminated  up to some number by using the commutativity property.
As we did earlier in the previous section, for each prime $l$ that divides $n$, we will consider $l$-part of the submodule and the symplectic group. Therefore we assume $n=l^e$ where $l$ is a prime number and let $W$ be a submodule on which there is an abelian action of a subgroup of the symplectic group.
\begin{equation*}
W={l^i\mathbb{Z}}/{l^e \mathbb{Z}}\times {l^{i+j_1}\mathbb{Z}}/{l^e\mathbb{Z}}\\
\times {l^{i+j_1+j_2}\mathbb{Z}}/{l^e \mathbb{Z}}\times...\times{l^{i+j_1+j_2+...j_{2m-1}}\mathbb{Z}}/{l^e\mathbb{Z}}.
\end{equation*}
For the action of the subgroup $\Gamma$ of the symplectic group to be well defined on the submodule, i.e $\Gamma W \subset W$, the entries of $\Gamma$ should satisfy some divisibility conditions. Hence the matrices in $\Gamma$ look like the following with those extra conditions on the entries due to sitting in the symplectic group.
$$ \left( \begin{array}{cccccc}
a_{11} & a_{12} & a_{13}& . & . & a_{12m}\\
l^{j_1}\bar{a}_{21} & a_{22} & a_{23}& . & . &a_{22m} \\
l^{j_1+j_2}\bar{a}_{31}& l^{j_2}\bar{a}_{32}&a_{33} & .& . &a_{32m} \\
. &.&.&.&.&.\\
. &.&.&.&.&.\\
l^{j_1+j_2+....j_{2m-1}}\bar{a}_{2m1} &l^{j_2+j_3+...j_{2m-1}} \bar{a}_{2m2} & l^{j_3+...j_{2m-1}} \bar{a}_{2m3}& .&.&a_{2m2m} \end{array} \right)$$
Let $k$ be the coulumn with $m$ "free" entries and $$J=\min \left \{e, \min ord_l (a_{k1}):  \left( \begin{array}{cccccc}
a_{11} & a_{12} & a_{13}& . & . & a_{12m} \\
a_{21} & a_{22} & a_{23}& . & . & a_{22m} \\
. &.&.&.&.&.\\
a_{k1} & a_{k2} & a_{k3}& . & . & a_{k2m}\\
. &.&.&.&.&.\\
a_{2m1} & a_{2m2} & a_{2m3}& . & . & a_{2m2m}
\end{array} \right )  \in\Gamma \right\}$$
Let $M$ be a matrix in the subgroup $\Gamma$ of the symplectic group with the $a_{k1}$'th entry divisible exactly by $l^J$, given as

$$M=\left( \begin{array}{cccccc}
a_{11} & a_{12} & a_{13}& . & . & a_{12m} \\
a_{21} & a_{22} & a_{23}& . & . & a_{22m} \\
. &.&.&.&.&.\\
l^{J}\bar{a}_{k1} & a_{k2} & a_{k3}& . & . & a_{k2m}\\
. &.&.&.&.&.\\
a_{2m1} & a_{2m2} & a_{2m3}& . & . & a_{2m2m}
\end{array} \right ) $$

And let $S$ be any matrix in $\Gamma$ represented as below.
$$S= \left( \begin{array}{cccccc}
a'_{11} & a'_{12} & a'_{13}& . & . & a'_{12m}\\
a'_{21} & a'_{22} & a'_{23}& . & . &a'_{2m} \\
. &.&.&.&.&.\\
l^{J}\bar{a}'_{k1}& {a}'_{k2}&a'_{k3} & .& . &a'_{km} \\
. &.&.&.&.&.\\
a_{2m1} &a'_{2m2} & a'_{2m3}& .&.&a'_{2m2m} \end{array} \right)$$

Since $\Gamma$ acts commutatively on $W$, we have

\begin{equation} \label{M}
 (MS-SM)  \left( \begin{array}{cccccc}
l^i& 0 & 0& . & . & 0 \\
0 & l^{i+j_1} & 0& . & . &0 \\
0& 0 &l^{i+j_1+j_2} & .& . &0 \\
. &.&.&.&.&.\\
. &.&.&.&.&.\\
0 & 0 & 0& .&.&l^{i+j_1+j_2+..j_{2m-1}}\end{array}\right) \equiv 0 \pmod {l^e}. 
\end{equation} 
For sake of simplicity in notation we will replace $\bar{a}'_{k1}$ with ${a}'_{k1}$ and $\bar{a}_{k1}$ with ${a}_{k1}$ that occurs in the representation of matrices $S$ and $M$. 

We will count the matrices $S$ in $\Gamma$. Consider the following relation. Given $S$ and fixed $M$, look at $X$' s that satisfy the following relation. By doing so we are decomposing the matrices $S$ in $\Gamma$ into components that we can use the implications of commutativity in counting.
\begin{equation}\label{forXsymp}
l^{i+J}[(a_{11}a'_{k1}-a'_{11}a_{k1})I+a'_{k1}S]=l^{i+J}[a_{k1}M+X]\,\pmod { l^e}
\end{equation}
$X$ has the form
\scriptsize
$$X=  \left( \begin{array}{cccccc}
0&a'_{k1}a_{12}-a_{k1}a'_{12}&a'_{k1}a_{13}-a_{k1} a'_{13}& . & . & a'_{k1}a_{12m}-a_{k1}a'_{12m} \\
a'_{k1}a_{21}-a_{k1}a'_{21} &a'_{k1}a_{11}-a_{k1}a'_{11}+a'_{k1}a_{22}-a_{k1}a'_{22}  & .& . & . &a'_{k1}a_{2m}-a_{k1}a'_{2m} \\
.&.&.&.&.&.\\
l^{J}(a'_{k1}a_{k2}-a_{k1}a'_{k2})&. & .& . &.&. \\
. &.&.&.&.&.\\
. &a'_{k1}a_{2m2}-a_{k1}a'_{2m2}) & 0& .&.&a'_{k1} a_{2m2m}- a_{k1} a'_{2m2m}\end{array}\right) \pmod{ l^{e-i-J}} $$
\normalsize  
We will use the implications of the assumption of commutativity to write some terms of $X$ in terms of others. For $1\leq n\leq 2m$. Let's assume $n\neq k$ and consider $l^{J}(a'_{k1}a_{kn}-a_{k1}a'_{kn})$ that occurs above in the representation of $X$.
Note that when the $n'th$ row of $(SM-MS)$ is multiplied with the first basis element of $W$ in \eqref{M} we get the first row first column of $\eqref{M}$.
\begin{equation}\label{nthrow} 
l^{i}((a'_{n1}a_{11}-a_{n1}a'_{11})(a'_{n2}a_{21}-a_{n2}a'_{21})+.....+\\
 l^{J}(a'_{nk}a_{k1}-a_{nk}a'_{k1})+....+(a'_{n2m}a_{m1}-a_{n2m}a'_{m1})).
\end{equation}
As the expression (\ref{nthrow}) is congruent to zero modulo $l^e$, we see that $(a'_{nk}a_{k1}-a_{nk}a'_{k1})$ is a combination of terms of $S$ outside the $k$'th column and $M$.
If $n=k$ the corresponding entry is $X_{kk}$ which is 
$$l^{J}(a'_{k1}a_{kk}-a_{k1}a'_{kk}+a_{11}a'_{k1}-a'_{11}a_{k1}).$$ 
Using the equivalence relation that comes from multiplying the $k$'th row of $(SM-MS)$ with the first basis element , we see that this entity can be written in terms of the entries in the columns beside the $k$'th one.  Hence all $a'_{nk}$ for $1 \leq n \leq 2m$ terms are determined by all the remaining columns' entries up to some multiple of $l^{(e-i-J)}$.
\\

From (\ref{forXsymp}) above, elements $S$ of $\Gamma$ can be decomposed into the sum of matrices as follows:
\begin{eqnarray}
S=xI+yM+X+ l^{e-i-J}Z \; \text{with}\;x,y\in {\mathbb{Z}}/{l^{e-i-J}\mathbb{Z}}, 
\,  \nonumber \\ Z\in M_{2m}({\mathbb{Z}}/{l^{i+J}\mathbb{Z}}) \text{ and } X \, \text{ as above}.
\end{eqnarray}

Note that the set of $Z$'s has dimension $2m^2+m$ since it's the $l^{e-i-J}$ part of the difference of two symplectic matrices, which comes from a set of dimension $2m^2+m$.
Also note that $x$, $y$ are determined by $a'_{11}$, $a'_{k1}$, $a_{k1}$, $a_{11}$ of $S$ and $M$.
\\
We have enough information to find a nice upper bound for $\|\Gamma \|$ which will be used to find a lower bound for $I(\Gamma)$.
$$
\begin{aligned}
|\Gamma| &\leq |\text{Set of S's}| \\&=
|{x}| \cdot|y|\cdot |X|\cdot |Z|\\&=
l^{(e-i-J)}l^{(e-i-J)}l^{(2m^2-2)(e-i-J)} l^{(2m^2+m)(i+J)}&=l^{(2m^2)e+m(i+J)}
\end{aligned}
$$
The order of the symplectic matrices of dimension $2m$ over the finite field $F$ where $|F|=q$ is known to be $|Sp(2m,F)|=q^{2m^2+m}\prod_{j=1}^{m}(1-q^{-2j})$ (see \cite{Dav}) and the order of $Sp(2m, Z/l^{e}Z)$  is derived as $|Sp(2m,Z/l^{e}Z)|=l^{e(2m^2+m)}\prod_{j=1}^{m}(1-l^{-2j})$(See Appendix).
This implies that,
\begin{equation}
\begin{aligned}\label{indexsub1}
I(\Gamma)&=\frac {\|Sp(2m,{\mathbb{Z}}/{l^e \mathbb{Z}})\|}{\| \Gamma \|}\\
&\geq \frac{l^{(2m^2+m)e}(1-l^{-2})(1-l^{-4})...(1-l^{-2m})}{l^{(2m^2)e+m(i+J)}}\\
&=l^{me-m(i+J)}(1-l^{-2})(1-l^{-4})....(1-l^{-2m}).
\end{aligned}
\end{equation}
On the other hand remember that the entries of matrices in $\Gamma$ are as the following,
\small
\[ \left( \begin{array}{ccccccc}
a_{11} & a_{12} & a_{13}& . & . &.& a_{12m}\\
l^{j_1}a_{21} & a_{22} & a_{23}& . & . &.&a_{22m} \\
l^{j_{(1)(2)}}a_{31}& l^{j_2}a_{32}&a_{33} & .& .&. &a_{32m} \\
. &.&.&.&.&.&.\\
l^{J}a_{k1} &.&.&.&.&.&.\\
. &.&.&.&.&.&.\\
l^{j_{(1)(2m-1)}} a_{2m1} &l^{j_{(2)(2m-1)}} a_{2m2} & l^{j_{(3)(2m-1)}} a_{2m3}& .&.&.&a_{2m2m}
 \end{array} \right).\]
\normalsize
Using the implications of the matrices being symplectic, we will try to see how many free entries (i.e entries taking all possible values in the ring across the symplectic group) we can get in the first column for matrices in $\Gamma$.  Using the representation of $S$ with block matrices with the property $A^T D-C^T B=I$ where $A^T D-C^T B$ is the matrix with entries $A_i D_j-C_i B_j$ with $1\leq i,j \leq m$ ,where $A_i,C_i$ are the $i$'th column of $A,C$; $D_j, B_j$ the $j$'th column of $D$ and $B$ respectively and the expressions $A_i D_j$, $C_i B_j$ are the inner product of two columns as vectors, we will figure out the range of the values of entries in the first column in the whole symplectic group. Assuming $B$ and $D$ are arbitrary fixed matrices we get $m$ linear relations on each column of $A,C$. $A,C$ has $2m$ elements in total at each column so at least $m$ of them takes all possible values in the ring of definition. Looking at the relation $A^T C=C^T A$, which is equivalent to $A_i C_j- A_j C_i=0$ for all $1\leq i,j \leq m$, we see that $m$ free entries in the first column of $A$ and $C$ should satisfy $m-1$ linear relations when the remaining columns of $A,C$  are fixed. 
The above relations on the first column $m_1$, as a row vector, of the symplectic matrix can be rewritten as:
$$m_1 \left ( \begin{array}{cc}
D & C'  \\
-B & -A' \end{array} \right)=v_1
$$ where $A',C'$ are $m\times(m-1)$ matrices obtained from $A$ and $C$ respectively by leaving out the first columns and $v_1$ is the row vector with $1$ in the first entry and $0$ for all remaining $2m-2$ entries. 
Since the matrix
$$\left ( \begin{array}{cc}
A & B  \\
C & D \end{array} \right)
$$ where $A,B,C,D$ are $m\times m$ matrices has $2m$ independent columns, the matrix
$$ \left ( \begin{array}{cc}
D & C'  \\
-B & -A' \end{array} \right)
$$ 
has $2m-1$ independent columns. Therefore the $(2m-1)$ relations on the first column of the symplectic matrix are independent and consistent. 
 
Hence we get one free entry in the first column.
Finding one free entry in the first column, all the other entries of the first column turn out to be a linear function of this entry. Without loss of generality it's possible to think any one of these $n$ entries free and express the others as a linear function of the free entry.  

Therefore we can assume that $a_{k1}$ of the symplectic matrices varies over the elements of the ring it's defined. Remember $J$ has been defined as the minimum of the valuations of $a_{k1}$ occurring in all matrices of $\Gamma$. The divisibility condition on the $k$th entry of the first column restrains the set of possible values of the free entry by $l^{J}$. Therefore we get the following bound on the index of $\Gamma$ in the symplectic group.

\begin{equation}
\begin{aligned}\label{indexsub2}
I(\Gamma) &\geq \frac{l^{(m^{2}+m)e}(1-l^{-2})(1-l^{-4})...(1-l^{-2m})}{l^{(m^2+m)e-J}}
\\&\geq l^{J}(1-l^{-2}) (1-l^{-4})...(1-l^{-2m}).
\end{aligned}
\end{equation}
Note that 
$\|W\|=l^{2me-2mi-(2m-1)j_1-c_2}$ where $c_2 \geq 0$.
The second power of (\ref{indexsub1}) multiplied with the $2m$'th power of (\ref{indexsub2})gives

$$I(\Gamma)^{2+2m}((1-l^{-2})...(1-l^{-2m}))^{-(2+2m)}\geq l^{2me-2mi}\geq |W|.$$
Since $\prod_l (1-l^{-2})^{-1}$ is $\zeta(2)$ and each $\prod_l (1-l^{-2t})^{-1}\leq \zeta(2)$ for $t\geq 1$ when we do not restrain $n$ in $Sp(2m,Z/nZ)$ to be only a prime power. 
We get
$$|W| \leq  I(\Gamma)^{2+2m}\zeta(2)^{(2+2m)}.$$
This ends the proof.

\section{Appendix}

Thanks to people at Mathoverflow site(S. Carnahan, J. Hahn, J. Humprehys, G. McNinch, C. Perez, A. Putman, J. Silverman, ) for showing me how to find the order of a symplectic group defined over a ring $Z/p^{k}Z$.  Here's the proof suggested:

Remember that

$$|Sp_{2n}(Z/pZ)|=p^{n^2}\prod ^{n}_{i=1}(p^{2i}-1)$$ 
and the reduction map  $ r: Sp_{2n}(Z/p^{k+1}Z)\mapsto Sp_{2n}(Z/p^{k}Z)$  is a surjective map.

To compute the order of the symplectic group it's enough to compute the size of the kernel of each reduction map. The elements in the kernel of the homomorphism $a:Sp_{2n}(Z/p^{k+1}Z)\mapsto Sp_{2n}(Z/p^{k}Z)$ have the form $I+p^{k} A$ for a symplectic matrix $A$. Also $I+p^k A$ satisfies the relation 
$$(I+p^k A)^{T}\Omega(I+p^k A)=\Omega \pmod {p^{k+1}},$$
where $\Omega$ is the matrix 
$$\left ( \begin{array}{cc}
0 & I_n  \\
-I_n & 0 \end{array} \right)
$$ with $I_n$, the $n\times n$ identity matrix.

The number of such $A$'s is $p^{2n^2+n}$. Therefore the order of the symplectic group over a finite ring $Z/p^{k}Z$ is given as, 
\begin{align}
|Sp_{2n}(Z/p^{k}Z)|=& p^{(2k-1)n^2+(k-1)n}\prod^{n}_{i=1}(p^{2i}-1)\\
=& p^{k(2n^2+n)}\prod^{n}_{i=1}(1-p^{-2i})
\end{align}




\end{document}